\documentclass[12pt]{paper}
\usepackage{color}
\usepackage{epsfig,cmmib57}
\usepackage{amssymb,amsmath,exscale}
\newcommand{\ZEP}{\epsilon}
\newcommand{\fn}{\phi_n'(0)}
\newcommand{\n}{\lambda_n}

\newcommand{\zsumU}{\sum _{n=1}^{+\ZIN}}
\newcommand{\ztz}{\tilde \zeta_n}

\newcommand\CM[1]{\vskip3mm\fbox{\parbox{5in}{#1}}\vskip3mm}
 
  \newcommand{\intx }{\int_0^x}

  \newcommand{\zzn}{\zeta_n}

\newcommand{\zg}{\gamma}

\newcommand{\intp}{\int_0^{\pi}}
\newcommand{\intT}{\int_0^T}
\newcommand{\intt}{\int_0^t}
\newcommand{\ints}{\int_0^s}
\newcommand{\intr}{\int_0^r}

\newcommand{\ZCD}{(\cdot)}
\newtheorem{Theorem}{Theorem}
\newtheorem{Corollary}[Theorem]{Corollary}
\newtheorem{Lemma}[Theorem]{Lemma}

\newtheorem{Remark}[Theorem]{Remark}
\newcommand{\ZOM}{\omega}
\newcommand{\zaa}{\alpha}
\newcommand{\ZR}{\rangle}
\newcommand{\ZL}{\langle}
\newcommand{\ZDE}{\delta}

\newcommand{\zdia}{~~\rule{1mm}{2mm}\par\medskip}

\newcommand{\ZLA}{\label}
\newcommand{\ZIN}{\infty}
\newcommand{\zProof}{{\noindent\bf\underbar{Proof}.}\ }

\newcommand{\ZBI}{\bibitem}
\newcommand{\ZD}{\;\mbox{\rm d}}
\newcommand{\zl}{\lambda}
\newcommand{\ZSI}{\sigma}

\author{
L. Pandolfi\thanks{Dipartimento di Scienze Matematiche, Politecnico di Torino, Corso Duca degli Abruzzi 24---10129 Torino, Italy (luciano.pandolfi@polito.it)}
}
\title{Boundary controllability and source reconstruction in a viscoelastic string under external traction\thanks{
This papers fits into the research program of the GNAMPA-INDAM and has been written in the framework of the   ``Groupement de Recherche en Contr\^ole des EDP entre   France et  Italie (CONEDP-CNRS)''.}}

\begin{document}
\maketitle
\pagestyle{myheadings}

\thispagestyle{plain}
 %\markright{L. Pandolfi}
 \markboth{L. Pandolfi }{Viscoelastic  string under external traction} 
\begin{abstract}
Treatises on vibrations devote large space to study the dynamical behavior of an elastic system subject to known external tractions. In fact, usually a ``system'' is not an isolated body but it is part of a chain of mechanisms which disturb the ``system'' for example due to the periodic rotation of shafts.
This kind of problem has been rarely studied in   control theory. In the specific case we shall study, the case of a viscoelastic string, the effect of such external action is on the horizontal component of the traction, and so it affects the coefficients of the corresponding wave type equation, which   will be time dependent. The usual methods used in controllability are not naturally adapted to this case. For example at first sight it might seem that moment methods can only be used in case of coefficients which are constant in time. Instead, we shall see that moment methods can be extended to study controllability of a viscoelastic string subject to external traction and in particular we shall study a controllability problem which is encountered in the solution of the inverse problem consisting in the identification of a distributed disturbance source.

\end{abstract}
\begin{keywords}    observability/controllability, integrodifferential system, moment problem, viscoelasticity      \end{keywords}
%\begin{AMS}Primary: 35Q93,45K05; Secondary:93B03 \end{AMS}

\section{Introduction}

Elastic and viscoelastic systems have been widely studied from the point of view of controllability (see for example~\cite{LAsTRI2,TUKSweiss,ZUAZUA} and references therein), the simpler ``canonical'' case being the string equation
\[ 
w_{tt}(\xi,t)=  \left (P c(\xi) w_{\xi }(\xi,t)\right )_{  \xi}
 \]
  The coefficients $ c(x)=1/\rho(x)  $ depends on the physical properties of the string  i.e. the density, and $ P $ is the orizontal traction in the string  which, for small oscillations and isolated systems, is constant in   time and space (while $ c $ may depend on $ x $ since the density may not be constant.)
  It is usually explicitly assumed in control problems that, a part the action of external controls which does not change the horizontal component of the stress,  the string is ``isolated'' i.e. that $ P $ does not depend on time.  Control problems for the string equations under external traction have been rarely studied, see~\cite{AB,AB-Bas,ABPand}. 
 In fact, in general the system is part of a chain of mechanisms which produce ``perturbations'' which influence the horizontal traction in the string, often in a known way, for example due to the rotation of shafts. For this reason, books   
 on vibrations devote large spaces to the study of the behaviour of elastic systems under known external tractions, see for example~\cite{Harto,Timo}.

In this paper we are going to study a  {\em viscoelastic\/} string subject to exterior known tractions. The string is also subject to a known signal $ \ZSI(t) $
 which enters trough an input operator $ b $ which is unknown. For example $ b $ may be the characteristic function of a certain set and our final goal will be the identification of the input operator $ b $ using boundary observations. So, the problem we are going to study boils down to the following equation
 
 \begin{equation}
\ZLA{eq:insECONDordine}
w_{tt} =P(t)\left (c(\xi)w_{\xi }\right) _\xi+\intt M(t-s) P(s) \left (c(\xi)w_{\xi }\right )_\xi(s)\ZD s +b\ZSI(t)
\end{equation}
where $ P(t) $ represents the effect of the variable external traction and $ w=w(\xi,t) $ with $ t>0  $ and $ \xi\in (0,\pi) $ (independent variables in $ w $ are not indicated unless needed for clarity).

\begin{Remark} 
{ \em
Writing the equation in the form~(\ref{eq:insECONDordine}) implies that the string is at rest for negative times. If not, a known additional terms will appear in the right hand side of the equation, which has no influence for the arguments below. So, we put it equal to zero.\zdia

}
\end{Remark}

Initial and boundary conditions are homogeneous (of Dirichlet type)
\begin{equation}\ZLA{eq:boundaryINItQUNAdisonoOmoge}
w(\xi,0)=0\,, \qquad w_t(\xi,0)=0\,,\qquad w(0,t)=w(\pi,t)=0\,.
 \end{equation}
The observation is the traction at $ x=0 $, i.e.

\[ 
 -P(t)\left [  c(\xi)w_\xi( 0, t)+\intt M(t-s)c(\xi)w_\xi( 0, s)\ZD s\right ]\,.
 \]

We shall assume that $ P(t) $, $ M(t) $ and $ c(x) $ are smooth  and that$ P(t) $ and $ c(x) $ are strictly positive (see below for the assumptions) so that from this ``real'' observation we can reconstruct the function
\begin{equation}
\ZLA{eq:osservMODIFICATA}
\eta(t)=w_\xi(0,t)
\end{equation}
which will be used in the reconstruction algorithm.

In the purely elastic case, this kind of inverse problem has been studied using control ideas in~\cite{Yamam} (see also~\cite{YamaGrasselli,PandIDENT}) where it is proved that this inverse problem (in the purely elastic case) depends on the solution of a boundary control problem. 
So, it is easily guessed that the control problem has to be studied first also in the case that the traction varies in time. This  is the second main subject  of this paper which is studied  first. 

{\em The problem we need to study first is the controllability of the deformation $ w(\cdot,T) $, at a certain time $ T $,\/} when $ b=0 $, the initial conditions are zero but the boundary conditions are
\begin{equation}
\ZLA{eq:boundPERcontroll}
  w(0,t)= \frac{f(t)}{c(0)P(t)}\,,\qquad w(\pi,t)=0 \,. 
\end{equation}
Here $ f(t) $ is a square integrable control (the denominator is introduced solely for convenience).  So, the control problem we shall study first, and which accupies the most part of this paper, is as follows: given any prescibed target $ W(\xi)\in L^2(0,\pi) $, we must prove the existence of a square integrable control $ f(t)\in L^2(0,T) $ such that $ w(\xi,T)=W(\xi) $.
Furthermore, we shall see that $ T $ does not depend on $ W $. 

The organization of the paper is as follows: the control problem is studied in section~\ref{sec:COntrol}   while the identification problem is studied in section~\ref{sec:identification}. Assumptions and preliminaries are in Section~\ref{sect:assump-references}.

\subsection{Further references}
Control problems for viscoelastic materials have already been studied under several different assumptions. We cite in particular~\cite{Kim2,ZANGPRIMOLAVORO} which seems to contain the most general results. The paper~\cite{Kim2} studied the case $ M=M(t,s) $ but it is explicitly stated that the arguments require constant density. The paper~\cite{ZANGPRIMOLAVORO} studied the case $ M=M(t,s,x) $ when the control is distributed (from which controllability results for boundary control systems with the control acting on the full boundary can be derived). These papers study also the multidimensional case, but the results does not seems to be ``constructive'' since they depend on compactness arguments (paper~\cite{Kim2}, which explicitly requires constant density) or Carleman estimates (paper~\cite{ZANGPRIMOLAVORO}). In this respect we cite also~\cite{PandAMO}, where controllability under boundary control is proved for multidemensional systems. The arguments are presented in the convolution case and constant density, but are easily extended to non constant density. The proof in this paper being based on compactness arguments is not constructive either. The paper~\cite{BarbuIANNELLI}, being based on an extension of D'Alambert formula (see last section) is more constructive in spirit, but the elastic operator is assumed to be time invariant.
The papers~\cite{PandIEOT,PandDCDS1,LorePANDsforza} instead study the control problem for a viscoelastic string, using moment methods, so that a representation formula for the control can be given (see in particular\cite[Section~4]{PandIEOT} and~\cite[Section~4]{LorePANDsforza}). Our goal here is the extension of these results to a viscoelastic string subject to nonconstant traction, and to show that the results can be used to solve a source identification problem. Source identification problems using control methods have been solved for elastic
 materials with constant traction first in~\cite{Yamam} (see also~\cite{YamaGrasselli}) and the method has been extended to materials with memory in~\cite{PandIDENT}. When the traction is not constant we need a different idea, presented in section~\ref{sec:identification}.

\section{\ZLA{sect:assump-references}Assumptions and preliminaries}
The assumptions in this paper are:
\begin{enumerate}
\item $ P(t) $ is continuous and strictly positive, $ P(t)>p_0>0 $ for every $ t\geq 0 $ while $  c(\xi) \in C^1(0,\pi) $ is also strictly positive: $  c(\xi) >  c_0  >0 $ for every $ \xi\in[0,\pi] $.
\item The kernel $ M(t)\in W_{\rm loc}^{2,2}(0,+\ZIN) $.  Hence its primitive $ N(t) $
\[ 
N(t)=1+\intt M(s)\ZD s
 \]
has three derivatives ($ N'''(t) $ is locally square integrable) and $ N(0)=1 $.
\end{enumerate}

  We introduce the following selfadjoint   operator $ A $:
  \begin{equation}\ZLA{eq:dFIOPeratorA}
{\rm dom}\, A= H^2(0,\pi)\cap H^1_0(0,\pi)\,,\qquad A\phi=
\left (c(\xi)\phi _{\xi }(\xi)\right )_\xi \,.
\end{equation}
It is well known that this operator has compact resolvent, and that it has a sequence of   eigenvalues $ \{-\zl_n^2\}_{n\geq 1} $, with the following asymptotic estimate~\cite[p.~173)]{Tricomi}
\begin{equation}
\ZLA{eq:stimeAUTOV}    \zl_n=n+\frac{H_n}{n}\,,\qquad |H_n|<M 
\end{equation}

A corresponding sequence of real normalized eigenfunctions $\phi_n(\xi)$ has the
following estimates, see~\cite[p.~176-177]{Tricomi}:
\begin{eqnarray}
&&\nonumber%\ZLA{eq:stimaPIUprecisaAutof}
\phi_n(\xi)=\sqrt{\frac{2}{\pi}} \left( \sin n\xi-\frac{1}{n}
K(\xi)\cos n\xi\right)+\frac{M_n(\xi)}{n^2}\,,\\
\ZLA{eq:stimaPIUprecisaAutof}&&=\sqrt{\frac{2}{\pi}}\left( \sin
n\xi +\frac{k_n(\xi)}{n}\right)\,,\quad \left\{ \begin{array}{ll} |K(\xi)|<M\,,&
|K'(\xi)|<M\,,\\ |k_n(\xi)|<M\,,& |M_n(\xi)|<M\,,
\end{array}\right.
\\
\ZLA{eq:stimeAUTOFDeriv}
 &&\phi'_n(\xi)=\sqrt{\frac{2}{\pi}}\left\{n\cos
n\xi+ \tilde k_n(\xi) \right\}\,,\qquad |\tilde k_n(\xi)|<M\,.
\end{eqnarray}

The previous formulas have been taken from~\cite{Tricomi} but this reference has $ n+1 $, with $ n\geq 0 $, where we wrote $ n $ (and we intend $ n\geq 1 $).
 Sharper estimates are
known, see~\cite[formula~(82)~p.~173 and~(85) p.~177]{Tricomi}  but we don't need this.

 \subsection{Riesz systems}
 
 The arguments of this paper are based on the theory of Riesz systems. We recall the definition: a sequence $ \{h_n\} $ in a separable Hilbert space $ H $ is a {\em Riesz basis\/} when it is the image of an orthonormal basis of $ H $ under a linear bounded and boundedly invertible transformation. The key property is that every $  h\in H$ ha a unique representation with respect to the Riesz basis $ \{h_n\} $,
 \[ 
 h=\sum \zaa_n h_n
  \]
  and $ \{\zaa_n\} \in l^2$. In fact, there are positive numbers $ m  $ and $ M $ such that
  \[ 
  m\|h\|^2\leq \sum |\zaa_n|^2\leq M\|h\|^2\,.
   \]
   If a sequence is a Riesz basis of its closed span (possibly not equal to $ H $) it is called a {\em Riesz sequence.\/}
   
   A property we shall need is that
  if $ \{h_n\} $ is a Riesz sequence then the series $ \sum \zaa_n h_n $ converges if and only if $ \{\zaa_n\}\in l^2 $, and the convergence is unconditional. Furthermore we need two perturbation results from~\cite[Ch.~1]{Young}. The first is a Paley-Wiener theorem, adapted to Hilbert spaces and orthonormal bases. It states that if
 \[ 
 \sum \|h_n-\ZEP_n\|^2<1
  \]
  ($ \{\ZEP_n\} $ orthonormal) 
  then $ \{h_n\} $ is a Riesz basis.  In fact we shall use the following corollary:
  \begin{Corollary}\ZLA{coro:estePaleyWIE}{\em
  Let $ \{\tilde \ZEP_n\} $ be a Riesz sequence and let the sequence $ \{h_n\} $ satisfy
  \begin{equation}
\ZLA{eq:PaleyWien}
 \sum \|h_n-\tilde \ZEP_n\|^2<+\ZIN
    \end{equation}
  then there exists a number $ N $ such that $ \{ h_n\} _{n>N} $ is a Riesz sequence too. This in particular implies that
\[%\begin{equation}
\ZLA{eq:PrimaCONDbari}
 \sum \zaa_n h_n
\]%\end{equation}
  converges in the norm of $ H $ if and only if $ \{\zaa_n\}\in l^2 $.
  }
   \end{Corollary}
   We stress that the sequence $ \{\tilde\ZEP_n\} $ in Corollary~\ref{coro:estePaleyWIE} need not be orthonormal.
 
Condition~(\ref{eq:PaleyWien}) does not imply that $ \{h_n\} $ is a Riesz sequence but 

\begin{Theorem}[Bari Theorem] {\em If the sequence $ \{h_n\} $ satisfies condition~(\ref{eq:PaleyWien}) and if furthermore
\begin{equation}
\ZLA{eq:deFiZOMindep}
\sum \zaa_n h_n=0\ \implies\ \{\zaa_n\}=0 
\end{equation}
then $ \{h_n\} $ is a Riesz sequence. 
}
\end{Theorem}

The additional condition~(\ref{eq:deFiZOMindep}) is called $ \ZOM $-independence.
Note that the convergence of the series in~(\ref{eq:deFiZOMindep}) is in $ H $ so that    it must be $ \{\zaa_n\}\in l^2 $ because $ \{h_n\}_{n>N} $ is a Riesz sequence.

\subsection{Representation of the solution}

Now we consider Eq.~(\ref{eq:insECONDordine}) with zero initial conditions but both the affine term and the boundary condition~(\ref{eq:boundPERcontroll}). We need a definition/representation formula for the solutions. We confine ourselves to a formula for $ w(x,t) $ and we ignore velocity, since  only $ w(\cdot,T) $ is needed for the identification problem. So, we prefer to write Eq.~(\ref{eq:insECONDordine})
 as a first order equation, integrating both the sides. We get
 \begin{equation}
\ZLA{eq:FORMAprimORDINE}
w_t=\intt N(t-s)P(s)\left ( c(\xi) w_{\xi }\right )_\xi(s)\ZD s+b(\xi) g(t)\,,\qquad g(t)=\intt \ZSI(s)\ZD S\,.
\end{equation}
%%%%%%%%%%%
 We multiply both sides of~(\ref{eq:FORMAprimORDINE}) with $ \phi_n(\xi)$ and we integrate on $ [0,\pi] $. We get
\begin{equation}
 \ZLA{eq:DIwN}
w_n'(t)=-\n^2\intt N(t-s)P(s) w_n(s)\ZD s 
 +\fn\intt N(t-s)f(s)\ZD s + b_n g(t) \\
 \end{equation}
 where
 \[
 w_n(t)=\intp w(\xi,t) \phi_n(\xi)\ZD x\,,\qquad b_n=\ZL b,\phi_n\ZR=\intp b(\xi)\phi_n(\xi)\ZD \xi \,.
\]

Let us fix   any time $ T>0 $. We want a representation formula for $ w_n(T) $. Let
\[ 
Q(t)=P(T-t) 
 \]
and let  $ z_n(t;T) $   solve
\[ 
z_n'(t;T)=-\n^2Q(t)\intt N(t-s)z_n(s;T)\ZD s\,,\qquad z_n(0)=1\,.
 \]
 It is simple to check the following representation formula for $ w_n(T) $:
 \begin{eqnarray}
\nonumber&&w_n(T)=\intT f(T-s)\left [
 \fn\ints N(s-r) z_n(r;T)\ZD r
 \right ]\ZD s\\
 \ZLA{eq:RappreFPRMUdiWn}&&+b_n \intT z_n(T-s;T)g(s)\ZD s\,. 
\end{eqnarray}
  so that we expect the following definition of the solution $ w(\xi,T) $:
 \begin{eqnarray} 
 \nonumber&&
  w(\xi,T)=\sum _{n=1}^{+\ZIN} \phi_n(\xi) \intT f(T-s)\left [
 \fn\ints N(s-r) z_n(r;T)\ZD r
 \right ]\ZD s\\
 \ZLA{eq:RAPPRESEseriediW}
 &&+
 \zsumU \phi_n(\xi) b_n \intt z_n(T-s;T)g(s)\ZD s\,.
   \end{eqnarray}
The following Lemma justifies formula~(\ref{eq:RAPPRESEseriediW}):
 \begin{Lemma}{\em 
\ZLA{eq:Lemma:regol-soluzione}
The series in~(\ref{eq:RAPPRESEseriediW}), as functions of $ T $, converge in $ C([0,K];L^2(0,\pi)) $ for every $ K\geq 0 $.
}
\end{Lemma}
The proof of this Lemma requires a preliminary study of the functions $ z_n(t,T) $ and can be found in the appendix.
%%%%%%%%%%K\geq 0
\section{\ZLA{sec:COntrol}The control problem}

In this section we study the control problem   needed to identify the source term $ b $; i.e. we consider Eq.~(\ref{eq:insECONDordine}) with $ b=0 $ and initial boundary conditions
\begin{equation}
\ZLA{eq:BOUNDARYcondiPeRcontrolPROBLEM}
w(\xi,0)=0\,,\quad w_t(\xi,0)=0\,,\qquad w(0,t)= \frac{f(t)}{c(0)P(t)}\,,\qquad w(\pi,t)=0  \,.
\end{equation}
 
As we already said, we want to identify a time $ T $ such that for every $ W(\xi)\in L^2(0,\pi) $ there exists a control $ f(t) $ such that the corresponding solution $ w(t) $ satisfies $ w(T)=w(\xi,T)=W(\xi ) $. 

It is convenient to perform several transformations first. In particular, as we are interested to the controllability solely of the deformation $ w(\cdot;T) $ and not of the velocity, we prefer to rewrite Eq.~(\ref{eq:insECONDordine}) in the form of a first order equation, but we note explicitly that the technique we use here can also be used to study the controllability of the pair $ (w(\cdot,T),w_t(\cdot,T)) $  as done in~\cite{LorePANDsforza} in the case that the external traction $ P $ and the density are constant. 

This section on controllability consists of three subsections: first we reduce the control problem to a moment problem. Then we  present the transformations which will simplify the problem uder study, in subsection~\ref{subs-CONTROLandtransform}. Then we prove that a certain sequence of functions is a Riesz sequence in $ L^2(0,S) $ wher $ S $ is a suitable number to be identified (in
 Theorem~\ref{teo:RieszSEQUENCEinTENpox}) and this will prove controllability, as seen in the 
subsection~\ref{subs-CONTROLandtransform}. The final result of controllability is as follows (compare~(\ref{eq:LaDefiS_0T0})):
\begin{Theorem}\ZLA{teo:ControllabINTEMpot}{\em 
Let $T_0 $ solve
\[ 
\int_0^{T_0} \sqrt{Q(s)}\ZD s =
\int_0^{T_0} \sqrt{P(T-s)}\ZD s
 =\pi 
 \]
and let $ T\geq T_0 $.
 Then, for every $ \xi\in L^2(0,\pi) $ there exists a boundary control $ f\in L^2(0,T) $ such that the corresponding solution $w(t) $ satisfies $w(T)=\xi $.
 }
\end{Theorem}
 
 Note that the number $ T_0 $ exists since
 \[ 
 \int_0^{T } \sqrt{P(T-s)}\ZD s\geq Tp_0\,.
  \]
  
The computation of $ w(t) $ for every $ t $ makes sense tanks to Lemma~\ref{eq:Lemma:regol-soluzione}.

 \begin{Remark}
{\em  In this Section~\ref{sec:COntrol} the final time $ T $ is fixed once and for all. So, instead of the notation $ z_n(t,T) $ we shall use the simpler notation $ z_n(t) $, i.e. we put $ z_n(t)=  z_n(t,T)$ ($ T $ fixed). We shall return to the complete notation in Section~\ref{sec:identification} and in the Appendix.\zdia

}
\end{Remark}

\subsection{\ZLA{subs-CONTROLandtransform}Reduction to a moment problem}

  We fix a certain time $ T $ and we consider the solution $ w(\xi,T) $, at this fixed time, when $ b=0 $. We have, from~(\ref{eq:RappreFPRMUdiWn}) and~(\ref{eq:RAPPRESEseriediW}):

 \begin{equation}\ZLA{eq:seriediW}
  w(\xi,T)=\sum _{n=1}^{+\ZIN} \phi_n(\xi) \intT f(T-s)\left [
 \fn\ints N(s-r) z_n(r)\ZD r
 \right ]\ZD s\,.
   \end{equation}
   As we noted, $ z_n(t) $  instead of $ z_n(t;T) $ since $ T $ is fixed.

   Equality~(\ref{eq:seriediW}) shows that $ w(\xi,T)=W(\xi ) $ if we can solve the moment problem
 \begin{equation}
\ZLA{probleMOMENTI}
\intT f(T-s)\left [
 \fn\ints N(s-r) z_n(r)\ZD r
 \right ]\ZD s=W_n  
\end{equation}
where
\[ 
 W_n=\intp \phi_n(\xi) W(\xi )\ZD \xi\,.
 \]
 
 Note that the transformation of $   L^2(0,\pi) \ni W\mapsto \{W_n\}$ is an isometric isomorphism of $ L^2(0,\pi) $ and $ l^2 $.

 Our  main result is that this moment problem is solvable and that the solution $ f $ of~(\ref{probleMOMENTI}) (which has minimal $ L^2 $-norm) depends continuously on $ \xi $, since we have:
 \begin{Theorem}\ZLA{teo:PRINCIPdieRIESZ} 
 {\em 
Let $ T=T_0 $  be as in Theorem~\ref{teo:ControllabINTEMpot} and $ T\geq T_0 $.

The sequence 
\begin{equation}
\ZLA{eq:sequeDAFARErieszSCRITTAconZ}
\left\{
 \fn\intt N(t-r) z_n(r)\ZD r
 \right \}
\end{equation}
 is a Riesz sequence in $ L^2(0,T) 
 $.\/}
\end{Theorem}

In order to prove the theorem, we need to perform several transformations.

\subsubsection{\ZLA{sect:trabNSf}Transformations}
Again we recall that we work with a fixed time $ T $ and $ 0\leq t\leq T $ so that   $ z_n(t;T) $  is simply denoted $ z_n(t) $ but it is important for the Appendix that we don't forget  the dependence on $ T $.

In this section, we transform the equation of $ z_n(t) $. 

First we introduce
\begin{equation}
\ZLA{eq:trasfH}
H(t)=-N'(0)t-\log Q(t)\,,\quad
\zzn(t)=e^{H(t)}z_n(t)=e^{-N'(0)t}\frac{1}{Q(t)}z_n(t)\,.
 \end{equation}
 We shall see below the reason for this   misterious looking transformation which does depend on $ T $ since $ Q(t)=P(T-t) $.
 
 The function $ \zzn(t) $ verifies
 \[ 
 \zzn'(t)-H'(t)\zzn(t)=-\n^2 e^{H(t)}Q(t)\intt N(t-s) e^{-H(s)}\zzn(s)\ZD s
  \]
 and so also
 \begin{equation} 
\ZLA{eq:diZZN}
\zzn''-H'(t)\zzn'+\left (\n^2 Q(t)-H''(t)\right )\zzn=G(t) 
\end{equation}
where

\begin{eqnarray}
%\zzn''-H'(t)\zzn'+\left (n^2 Q(t)-H''(t)\right )\zzn=G(t) 
\nonumber&&
G(t)=-\n^2\left [
H'(t)Q(t)e^{H(t)}+e^{H(t)}Q'(t)
\right ]\intt N(t-s)e^{-H(s)}\zzn(s)\ZD s\\
\ZLA{eq:diGRECOzetaN}&&
-\n^2e^{H(t)}Q(t)\intt N'(t-s) e^{-H(s)}\zzn(s)\ZD s\,.
\end{eqnarray}

Now we apply Liouville transformation to Eq.~(\ref{eq:diZZN}), which requires the introduction of two auxiliary functions, which will be determined later on (see~\cite[p.~163]{Tricomi}):
 we introduce the function $ \ztz(t) $ given by
 \begin{equation}
 \ZLA{eq:LadefinIZdiZETAtilde}
a(t) \ztz(t)=\zzn(t)
  \end{equation}
 (here $ a(t) $ is the first auxiliary function). Then we have
 \begin{eqnarray*}
&&a(t)\ztz''(t)+\left [2a'(t)-H'(t)a(t)\right ] \ztz'(t)\\
&&+\left [
\left (\n^2Q(t)-H''(t)\right )a(t)-H'(t)a'(t)+a''(t)
\right ]\ztz(t)=G(t)\,.
\end{eqnarray*}
The second trasformation is a transformation of the time variable: we introduce
\[ 
x=L(t)
 \]
 where $ L(t) $ is still unspecified (but we shall see that we can choose $ L(t) $ strictly increasing and its inverse transformation will be denoted $ M(t) $.)  We introduce $ Y_n(x) $ defined by
 \[ 
\ztz(t)=Y_n(L(t))\,.
  \]
 Then we have the following equation for $ Y_n(x) $:\
 \begin{eqnarray*}
&&a(t)\left [L'(t)\right] ^2 Y''(L(t))+\left [ a(t)L''(t)+\left (2a'(t)-H'(t)a(t)\right )L'(t)\right ] Y'(L(t))\\
&&+ \left [ \n^2Q(t)a(t)+\left ( a''(t)-H''(t)a(t)-H'(t)a'(t)\right )\right ]Y(L(t))=G(t)\,.
\end{eqnarray*}
We must be precise on this point: here and below notations like $ Y'( L(t))$ will denote the derivative of $ Y(x) $ computed for $ x=L(t) $.

Now we relate the functions $ L(t) $ and $ a(t) $ so to have the coefficient of $ Y'(L(t)) $ equal to zero, i.e. we impose
\[
\left(L'(t)\right )' =-L'(t)\left [
2\frac{a'(t)}{a(t)}-H'(t)
\right ] 
\]
 (which is legitimate since we shall choose $ a(t)\neq 0 $.) In order to satisfy this condition we choose
 \[
 L'(t)=\frac{e^{H(t)}}{a^2(t)}\,.
   \]
   Note that $ L(t) $ turns out to be strictly  increasing. Furthermore, we choose
   \[ 
   L(t)=\intt \frac{e^{H(s)}}{a^2(s)}\ZD s
    \]
    so that we have also $ L(0)=0 $ and   $ L(t)>0 $ for $ t>0 $.
    
    Once $ L(t) $ has been chosen as above, we see that the equation of $ Y_n(x) $ is
    \begin{eqnarray*}
&& \frac{1}{a^3(t)}e^{2H(t)}Y_n''(L(t)) +\left [
\n^2Q(t)+\left ( a''(t)-H''(t)a(t)-H'(t)a'(t)\right )
\right ]Y_n(t)\\
&&=G(t)\,.
\end{eqnarray*}
Finally we choose 
\[ 
a(t)=e^{H(t)/2}\frac{1}{\sqrt[4]{Q(t)}}
 \]
 (note: it is strictly positive) and we get the final form of the equation for $ Y_n(t) $: 
 \begin{equation}
\ZLA{eqYnFORMAfinale}Y''(L(t))+\left [\n^2 +\tilde V(t)\right ] Y(L(t))=\frac{e^{-H(t)/2}}{Q(t)^{3/4}}G(t)
\end{equation} 
where
\[ 
\tilde V(t)=\frac{e^{-H(t)/2}}{Q(t)^{3/4}}\left [ a''(t)-H''(t)a(t)-H'(t)a'(t)\right ]\,.
 \]
The definition of $ a(t) $ gives
\begin{equation}
\ZLA{eq:ladefinizionediLprimo}
L'(t)=\sqrt {Q(t)}\,.
 \end{equation}
Let $ M(x) $ be the inverse function of $ L(t) $, defined on the interval $ [0,L(T)] =[0,S]$ and 
\[ 
V(x)=\tilde V(M(x))\,.
 \]
 Using this and the explicit form of $ G(t) $ in~(\ref{eq:diGRECOzetaN})  we see that $ Y_n(x) $ satisfies the following integro-differential equation
 on the interval
 \[ 
 0\leq x\leq S=L(T)\,: 
  \]
 \begin{eqnarray*}
&&Y''(x)+\left (\n^2+V(x)\right )Y_n(x)\\
&&=-\n^2e^{H(M(x))/2}\left [ H'(M(x)) Q^{1/4}(M(x))\right.\\
&&\left.+Q'(M(x))Q^{-3/4}(M(x))\right ] \int_0^{M(x)}N(M(x)-s)e^{-H(s)}\zzn(s)\ZD s \\
&&-\n^2 e^{H(M(x))/2}Q^{1/4}(M(x))\int_0^{M(x)} N'(M(x)-s) e^{-H(s)}\zzn(s)\ZD s\,. 
\end{eqnarray*}
We make the substitution $ s=M(r) $ (we recall that $ M(r) $ is the inverse function of $ L(s) $), we use $ 0\leq r\leq x $ when $ 0\leq s\leq M(x) $ and
\[ 
\frac{\ZD}{\ZD r}M(r)=\frac{1}{L'(s)}=\frac{1}{\sqrt{Q(s)}}\qquad \mbox{where $ s=M(r) $.}
 \]
 Then,   we have the following equation for $ Y_n=Y_n(x) $:
 \begin{equation}
 \ZLA{eq:FinalePERyNcompressa}
 Y_n''+V(x)Y_n+\n^2Y_n=-\n^2\intx A(x,r)Y_n(r)\ZD r 
 \end{equation}
 where 
 \begin{align}
 \nonumber
 A(x,r)=e^{H(M(x))/2}Q^{1/4}(M(x))\left \{
 \left [
 H'(M(x))\right.\right.\\
\nonumber 
 \left.  +Q'(M(x))Q^{-1}(M(x))
 \right ]N(M(x)-M(r)) %\\
 %\left.
 \\
 \ZLA{eq:FormaDiA(xr)}
 \left.
 +N'(M(x)-M(r)) 
 \right \}e^{-H(M(r))/2}Q^{-3/4}(M(r))\,.
 \end{align}

 The values $ Y_n(0)=y_0 $ and $ Y'_n(0)=y_1 $ are easily computed:
 \begin{eqnarray*}
&& y_0=Y_n(0)=e^{H(0)/2}\sqrt[4]{Q(0)}=Q (0)^{-1/4}=P (T)^{-1/4}\,,\\
&&y_1=Y_n'(0)=-\frac{1}{4}Q(0)^{-7/4}\left \{
2N'(0)Q(0)+Q'(0)
\right \}\\
&&=\frac{1}{4}P(T)^{-7/4}\left \{
P'(T)-2N'(0) P(T)
\right \}\,.
\end{eqnarray*}

 In the following, we don't need these explicit expressions, but we shall need the following facts:
 \begin{itemize}
\item the initial conditions $ y_0 $ and $ y_1 $ do not depend on $ n $;
\item the initial condition $ y_0  $ is strictly positive (we shall need that it is different from zero.)
\item the initial conditions $ y_0 $ and $ y_n $ do depend on $ T $, equivalently on $S=L(T)  $. 
 
 In the appendix we shall use the fact that when $ T\in[0,K] $ then there exists a number $ M $ (which depends on $ K $) such that $ |y_0(T)|<M $, $ |y_1(T)|<M $ for every $ T\in[0,K] $.
\end{itemize}
  We introduce the notation
  \begin{equation}
\ZLA{eq:definDIgN(x)}
g_n(x)=y_0\cos \n x+\frac{1}{\n}y_1\sin \n x\,.
\end{equation}
 
Using~(\ref{eq:definDIgN(x)}) and~(\ref{eq:FinalePERyNcompressa}) we can write the following Volterra integral equation for $ Y_n(x) $:
\begin{eqnarray*}
\nonumber
&&Y_n(x)=g_n(x) -\frac{1}{\n}\intx \sin \n(x-s) V(s) Y_n(s)\ZD s
\\
&&\nonumber
-\n\intx \sin \n(x-s)\ints A(s,\nu)Y_n(\nu)\ZD \nu\,\ZD s\\
&& \nonumber = g_n(x)
-\left \{
 \intx \left [A(x,s) +\frac{1}{\n} \sin \n(x-s) V(s)\right ] Y_n(s)\ZD s\right.\\
&&\nonumber   \left. + \intx \cos \n(x-s)A(s,s) Y_n(s)\ZD s\right.\\
&&\left.+\intx \cos \n(x-s)\ints A_{,1}(s,\nu)Y_n(\nu)\ZD\nu\,\ZD s\right\}\,.
\end{eqnarray*} 

In this and following formulas, we use the comma notation for the derivative. Hence,  $ A_{,1}(x,s) $ is the derivative of $ A(x,s) $ respect to the first variable.

When $ x=r $ the brace in~(\ref{eq:FormaDiA(xr)}) is  
\begin{eqnarray*}
 \left \{ H'(M(r))+Q'(M(r))Q^{-1}(M(r)) +N'(0)\right \}=0
\end{eqnarray*}
thanks to the choice~(\ref{eq:trasfH}) for $ H(t) $. 
The reason for the choice of the exponential $ e^{H(t)} $ is precisely this:
\[ 
A(s,s)=0\,.
 \]

So, elaborating further the integral equation for $ Y_n(x) $ we have
\begin{eqnarray}
\nonumber&& Y_n(x)=g_n(x)-\intx A_n(x,s)Y_n(s)\ZD s \\
\nonumber
%\ZLA{eq:PrimRApprREYn} 
&&+\frac{1}{\n}\intx\sin \n(x-s)\left [\ints A_{,11}(s,\nu)Y_n(\nu)\ZD\nu   \right ]\ZD s \\
%%%%%%%5
\nonumber&&=g_n(x)-\intx A_n(x,s)Y_n(s)\ZD s \\
\ZLA{eq:PrimRApprREYn} &&+\frac{1}{\n}\intx Y_n(\nu)\left [
\int_0^{x-\nu}  A _{,11}(x-r,\nu) \sin\n r\ZD r
\right ]\ZD\nu\,.
%%%%%%%
\end{eqnarray}
where   
\[
A_n(x,s)=\left [A(x,s) +\frac{1}{\n}\left ( V(s)-A_{,1}(s,s)\right ) \sin \n(x-s) \right ] \,.
\]
Note that 
\[ A_n(x,x)=0\,.
\]

\subsection{Some estimates}
The estimates in the following Lemma~\ref{lemmaVicinaNZAdiYn}  are used also in the Appendix, where we need to keep track of the 
dependence on   $ T\in[0,K] $, hence on $ S$ in the corresponding interval $[0,\Xi] $. So, in this lemma dependence on $ S $ is explicitly indicated. For this, we recall that the function $ Y_n(x)=Y_n(x;S )$ is defined for $ 0\leq x\leq S\leq \Xi $. Furthermore we recall that also $ y_0 $ and $ y_1  $   depend on $ S $: $ y_0=y_0(S) $  $ y_1=y_1(S) $. 
Using the asymptotic formulas~(\ref{eq:stimaPIUprecisaAutof}) and Gronwall inequality we get:
\begin{Lemma}
{\em 
\ZLA{lemmaVicinaNZAdiYn} 
Let $ \Xi>0 $. There exists a number $ M=M_\Xi $ such that for every $ x\in [0,S]\subseteq [0,\Xi] $ the following inequalities hold:
\begin{equation}
\ZLA{eq:VicinaNZAdiYn}
| Y_n(x;S)|<M\,,\qquad \left | Y_n(x;S)-y_0(S)\cos n x\right |\leq \frac{M}{n}\,. 
\end{equation}
}
\end{Lemma}
The proof is in the Appendix.

Using Corollary~\ref{coro:estePaleyWIE} we get:
\begin{Theorem}{\em 
Let $ S\geq \pi $ be fixed. There exists $ N>0 $ such that the sequence $ \{Y_n(x)\} _{n\geq N} $ is a Riesz sequence in $ L^2(0,S) $.
}
\end{Theorem}
This suggests that we can use Bari Theorem, as in~\cite{PandDCDS1}, in order to prove that $ \{Y_n(x)\} _{n\geq 1}  $ is a Riesz sequence. This is an intermediate step we shall need below, but we note that the sequence to be studied is~(\ref{eq:sequeDAFARErieszSCRITTAconZ}) which, written in terms of $ Y_n(x) $ gives the following sequence of functions of the variable $ t $:
\[ 
\fn\int_0^{L(t)}N(t-M(s)) e^{N'(0) M(s)/2}
\frac{1}{ \sqrt[4]{Q(M(s))} }Y_n(s)\ZD s\,.
 \] 
The transformation from $ \psi(t)\in L^2(0,T) $ to $( {\cal L}\psi)(x) =\psi\left (M(x)\right ) $ in $ L^2(0,S) $ (where $ T=M(S) $) defines a bounded and boundedly invertible transformation of $ L^2(0,T) $ onto $ L^2(0,S) $ so that 
  we can equivalently study the following sequence in $ L^2(0,S) $:
\begin{equation}
\ZLA{eq:laSUCCESSconCxs}
\left \{\fn\intx C(x,s) Y_n(s)\ZD s\right \}    
\end{equation}
where
\begin{equation}
\ZLA{eq:definCxs}
C(x,s)=
 N(M(x)-M(s)) e^{N'(0) M(s)/2} Q^{-1/4}(M(s)) 
 \end{equation}
and we must identify a value of $ S $ such that the sequence in~(\ref{eq:laSUCCESSconCxs}) is Riesz in $ L^2(0,S) $. 
We shall prove:
\begin{Theorem}\ZLA{teo:RieszSEQUENCEinTENpox}{\em 
The sequence~(\ref{eq:laSUCCESSconCxs}) is a Riesz sequence in $ L^2(0,S) $ for every $S\geq \pi$\,.
}
\end{Theorem}

Using~(\ref{eq:ladefinizionediLprimo}) and $ L'(0)=0 $ we see that the sequence~(\ref{eq:sequeDAFARErieszSCRITTAconZ}) is a Riesz sequence in $ L^2(0,T) $ for every $ T\geq M(\pi) $, i.e. for every $ T>T_0 $ where $ T_0 $ solve
\begin{equation}
\ZLA{eq:LaDefiS_0T0}
\int_0^{T_0} \sqrt{Q(s)}\ZD s=\pi\,.
 \end{equation}
  This is the statement in Theorem~\ref{teo:ControllabINTEMpot}.

Thanks to the fact that $ C(x,x) $ is strictly positive, we can equivalently study the  sequence
\begin{equation}
\ZLA{eq:sequeDAFARErieszSCRITTAconY}
\left \{\fn\intx B(x,s) Y_n(s)\ZD s\right \}\,, \qquad B(x,s)=
\frac{C(x,s)}{C(x,x)} \,.
\end{equation}
The computational advantage of this transformation is that now $ B(x,x)=1 $.

\subsubsection{The proof of Theorem~\ref{teo:RieszSEQUENCEinTENpox}}

In this proof we work on a fixed interval $ [0,S] $ (hence in a fixed interval $ [0,T] $ so that we don't need to take track of the dependence of   $ Y_n(x) $ on $ S $. We shall need   asymptotic estimates for several different sequences. In order to simplify the notations, $ M $ will denote  a number which does not depend on the index of the sequences  but which in general depends on the interval $ [0,S] $ we are working with. Analogously, $ \{M_n(x)\} $ will denote a  sequence of functions which is bounded on $ [0,S] $. The bound $ M $ can depend on $ S $.
      These constants and sequences will not be the same at every occurrence, without any possibility of confusion.

 We introduce 
\begin{equation}
\ZLA{eq:DefinizionediK}
K(s)=V(s)- A_{,1}(s,s)
\end{equation}
so that
$ A_n(x,s) $ in formula~(\ref{eq:PrimRApprREYn}) is
\[ 
A_n(x,s)=A(x,s)+\frac{1}{\n}K(s)\sin \n(x-s)\,.
 \]
Furthermore, we introduce the notation  $  Z_n(x)  $:   in~(\ref{eq:sequeDAFARErieszSCRITTAconY}):
\begin{equation}
\ZLA{eq:defiZn}
Z_n(x)=\fn\intx B(x,s) Y_n(s)\ZD s\,.
\end{equation}

We shall proceed in parallel with the sequences $ \{ Y_n(x)\} $ and $ \{ Z_n(x)\} $ in order to prove the following result  which, as we noted, implies Theorem~\ref{teo:PRINCIPdieRIESZ}.

\begin{Theorem}\ZLA{teo:Rieszbothsequences}{\em 
Let $ S\geq \pi $. Both the sequences $\{ Y_n(x)\} $ and $ \{Z_n(x)\} $ are Riesz sequences in $ L^2(0,S) $.
}
\end{Theorem}
The proof require several steps.

\subparagraph{Step 1: linear independence.}

We first prove that if $ \{Z_n(x)\} $ is linearly dependent then $ \{Y_n(x)\} $ is linearly dependent too; and then we prove linear independence of $ \{Y_n(x)\} $. 

If 
\[ 
\sum _{n=1} ^{K}\zaa_n Z_n(x)=0\quad {\rm  i.e.} \quad \intx B(x,s)\left [
\sum _{n=1} ^{K}\fn \zaa_n Y_n(s)
\right ]\ZD s=0
 \]
then we have  
\[
\sum _{n=1} ^{K}\fn \zaa_n Y_n(s)=0 
\]
 because $ B(x,s) $ is smooth and $ B(x,x)=1$.

Now we prove that
$ \{ Y_n(x)\} $ is linearly independent. We prceed by contradiction: let   $ K $ be  the first index for which
\begin{equation}
\ZLA{eq:LinearDEPENdiYn}
\sum _{n=1} ^{K}  \zaa_n Y_n(s)=0 \,.
\end{equation}  
Note that $ K\geq 2 $ since $ Y_1(x)\neq 0 $.

Using equality~(\ref{eq:FinalePERyNcompressa}) we get
\begin{eqnarray*}
0=\sum _{n=1} ^{K}\zaa_n Y''_n(x)=-\sum _{n=1} ^{K}\n^2\zaa_n Y_n-\intx A(x,r) \left (
\sum _{n=1} ^{K}\n^2\zaa_n Y_n(r)
\right )\ZD r
\end{eqnarray*}
So we have also 
\[ 
\sum _{n=1} ^{K}\n^2\zaa_n Y_n(x)=0\,.
 \]
 This equality and~(\ref{eq:LinearDEPENdiYn}) contradict the definition of $ K $ since they give:
 \[ 
 \sum _{n=1} ^{K-1}(\zl_K^2-\n^2)\zaa_n Y_n(x)=0\,.
  \] 

   \subparagraph{Step~2: a new expression for $ Z_n(x) $.} 
   
   We integrate by parts in~(\ref{eq:PrimRApprREYn}) as follows (the definition of $ K(s) $ is in~(\ref{eq:DefinizionediK})):

\begin{eqnarray*}
&&
   Y_n(x)=g_n(x) -\intx A(x,s)Y_n(s)\ZD s   -\frac{1}{\n}\intx \sin \n(x-s) K(s)Y_n(s)\ZD s\\
   &&
   -\frac{1}{\n^2}\intx\left [ \int_0^{x-\nu}A _{,11}(x-s,\nu) \frac{\ZD}{\ZD s}\cos \n s\ZD s\right ]Y_n(\nu)\ZD\nu\\
&&   = g_n(x)-\intx A(x,s)Y_n(s)\ZD s
    -\frac{1}{\n}\intx \sin \n(x-s) K(s)Y_n(s)\ZD s\\
    &&+\frac{1}{n^2}M_n(x)\,.
 \end{eqnarray*}

   Let $ R(x,s) $ be the resolvent kernel of $ -A(x,s) $.  Then we have (recall $ \n\sim n $)
   \begin{eqnarray}
\nonumber && Y_n(x)=\left [       g_n(x)
-\frac{1}{\n} \intx \sin \n(x-s) K(s)Y_n(s)\ZD s
\right ]\\
\nonumber&&+\intx R(x,s)\left [
g_n(s)-\frac{1}{\n}\ints K(\nu)\sin \n(s-\nu)Y_n(\nu)\ZD\nu
\right ]\ZD s\\
\ZLA{eqExprediYnCONrisolv}&&+\frac{1}{n^2}M_n(x)\,.
\end{eqnarray}
Using  smoothness of $  R(x,s) $ and $ R(x,x)=0 $ (a consequence of $ A(x,x)=0 $), we can integrate by parts twice and we see:
 
\[ 
\intx R(x,s)g_n(s)\ZD s=\frac{M_n(x)}{n^2}\,.
 \]
Using Lemma~\ref{lemmaVicinaNZAdiYn}, we see that the same holds also for the double integral (exchange the order of integration and integrate by parts once) so that
 
\[ 
Y_n(x)=\left [       g_n(x)
-\frac{1}{\n} \intx \sin \n(x-s) K(s)Y_n(s)\ZD s
\right ]+\frac{1}{n^2}M_n(x)\,.
 \]

    We replace this expression in~(\ref{eq:defiZn}) and we see that:
    \begin{eqnarray}
\nonumber &&Z_n(x)=\fn y_0\intx B(x,s) \cos \n s\ZD s+y_1\frac{\fn}{\n}\intx B(x,s)\sin \n s\ZD s\\
\nonumber&&-\frac{\fn}{\n}\intx \left [\int _{\nu}^x B(x,s)\sin \n(s-\nu)\ZD s\right ]K(\nu) Y_n(\nu)\ZD\nu\\
&& 
\ZLA{eq:PrimINTEPartiB}
+\frac{1}{\n}\intx B(x,s)M_n(s)\ZD s  \,.
\end{eqnarray}
  
Integrating by parts,  we see that on every interval $ [0,S] $ the functions $ Z_n(x)  $ are sum of a term of the order $ 1/n $ plus:
\begin{equation}
\ZLA{eq:intePARTiinB}
\fn y_0\intx B(x,s) \cos \n s\ZD s=y_0\frac{\fn}{\n}\left [\sin \n x- \intx B _{,2}(x,s)\sin \n s\ZD s\right ]\,.
\end{equation}
A further integration by parts  of the last integral 
and the estimates~(\ref{eq:stimeAUTOV}) and~(\ref{eq:stimeAUTOFDeriv})
show the following result:
\begin{Lemma}\ZLA{eq:STIMAbase-perZn}{\em 
Let $ S>0 $. There exists $ M $ such that for every $ n $ and every $ x\in[0,S] $ we have
\begin{equation}
\ZLA{eq:DISEQperZ}
|Z_n(x)|<M\,,\qquad |Z_n(x)-y_0\sqrt{\frac{2}{\pi}} \sin n x |<\frac{M}{n}\,.
 \end{equation}
 }
\end{Lemma}

 \subparagraph{Step 3: we prove that $ \ZOM $-independence of $ \{Y_n(x)\} $ on   $L^2(0,S) $, $S\geq \pi$, implies that of $ \{Z_n(x)\} $.}
 
 We proceed by contradiction. Let $\{\zaa_n\}$ be a sequence such that the following equality holds in $ L^2(0,S) $:
\begin{equation}
\ZLA{eq:SeZnEdipende}
\zsumU \zaa_n Z_n(x)=0\,.
\end{equation}
 The fact that $ \{ Z_n(x)\} _{n>N} $ is a Riesz sequence in $ L^2(0,S) $ implies that $ \{\zaa_n\}\in l^2 $. In fact, we have more:
 \begin{Lemma}  \ZLA{lemma:PrimoINappendice}{\em 
If~(\ref{eq:SeZnEdipende}) holds  in $ L^2(0,S) $ (any $ S>0 $) then the series 
\begin{equation} \ZLA{eq:LaserieMOLTOpdisinn x}
\sum \fn\zaa_n\sin n x
 \end{equation}
 converges in $ L^2(0,S) $.
So, 
If~(\ref{eq:SeZnEdipende}) holds and if $ S\geq \pi $ then there exists $ \{\zg_n\}\in l^2 $ such that $$ \zaa_n=\frac{\zg_n}{\fn}\,. $$
}
\end{Lemma}   
\zProof 
The second statement follows from the first one, since the sequence $ \{  \sin\n x\}_{n>N} $
 is a Riesz sequence\footnote{here in fact $ N=1 $ but we don't need this.} in $ L^2(0,S) $ when $ S\geq \pi $ so that convergence of~(\ref{eq:LaserieMOLTOpdisinn x}) implies $ \{\fn\zaa_n\}\in l^2 $. So, we prove that equality~(\ref{eq:SeZnEdipende}) implies convergence of the series~(\ref{eq:LaserieMOLTOpdisinn x}).
 
We prove the first statement. Using~(\ref{eq:PrimINTEPartiB}), equality~(\ref{eq:SeZnEdipende}) can be written as
\begin{eqnarray*}
&& \zsumU \frac{\fn}{\n}\zaa_n\intx \left [
\int _{\nu}^x B(x,s)\sin \n(s-\nu)\ZD s
\right ]K(\nu)Y_n(\nu)\ZD\nu\\
&&-\zsumU \frac{\zaa_n}{\n}\intx B(x,s) M_n(s)\ZD s
-y_1 \intx B(x,s)\left (
\zsumU \frac{\fn}{\n} \zaa_n \sin \n s\ZD s\right ) 
  \\
&&= y_0\left (
\zsumU  \fn\zaa_n \intx B(x,s)\cos \n s\ZD s
\right )
\end{eqnarray*}
(these equalities technically have to be intended as finite sums, and the equality holds in the limit if it happens that the series converges. But, every series a part possibly the one on the right hand side   is clearly convergent, so that the series on the right side has to converge too).
Using~(\ref{eq:intePARTiinB}) we write the series on the right hand side as   
 \[ 
 y_0\left (
 \zsumU \zaa_n\frac{\fn}{\n}\sin \n x
 \right )-y_0\intx B_{,2}(x,s)\left (
 \zsumU \zaa_n\frac{\fn}{\n}\sin \n s\ZD s
 \right )\ZD s
  \]
We deduce for this  that every term, a part possibly the first series above, can be differentiated in $ L^2 $ and that the differentiated series converge. Hence  also the first series can be differentiated termwise and
\[ 
y_0\left (
 \zsumU \zaa_n\fn \cos \n x
 \right )
 \]
converges to a square integrable function, so that $ \{\fn\zaa_n\}\in l^2 $ since $ \{\cos \n x\}_{n>N} $ is a Riesz sequence in $ L^2(0,S) $ when $ S\geq\pi $.\zdia
  
  Consequently we have
  \[ 
  0=\zsumU \frac{\zg_n}{\fn} Z_n(x)=\intx B(x,s) \left [\zsumU \zg_n Y_n(s)\right ]\ZD s\,.
   \]
   The fact that $ B(x,s) $ is smooth and $ B(x,x)=1 $ easily implies
\[
   \zsumU \zg_n Y_n(s)=0\,.
\]

   Arguing by contradiction, we see that if $ \{Y_n(x)\} $  is $ \ZOM $-independent, the same holds of $ \{Z_n(x)\} $.    
 So, in order to prove that $ \{Z_n(x)\} $ is $ \ZOM $-independent, which is our final goal, it is sufficient to prove $\ZOM  $-independence of $ \{ Y_n(x)\} $. This we do in the last step.
  
 \subparagraph{Step 4: the sequence $ \{ Y_n(x)\} $ is $ \ZOM $-independent.}
 
Let the sequence $ \{\zg_n\} $ satisfy
\begin{equation}
\ZLA{eq:PseKAYnFOSedipende}
\zsumU \zg_nY_n(x)=0
 \end{equation}
 in $ L^2(0,S) $ (so that $ \{\zg_n\}\in l^2 $). Using~(\ref{eq:PrimRApprREYn}) we get
 \begin{eqnarray*}
&& \sum \zg_n g_n(x)=
 \sum \zg_n \intx A_n(x,s)Y_n(s)\ZD s
 \\
 &&-\sum \frac{\zg_n}{\n^2}\left \{\intx 
 A_{,11}(x,\nu)Y_n(\nu)\ZD\nu \right.\\
 &&\left. +\intx\cos \n(x-s)\left [ A_{,11}(s,s)Y_n(s)  +\int_{x-s}^x A_{,111}(s,\nu) Y_n(\nu)\ZD\nu\right ]\ZD s\right \}\,. 
 \end{eqnarray*}
 We prove  that we can termwise compute the derivative of the series on the left hand side. As above, it is sufficient that  we prove that the series on the right side can be differentiated termwise.  This is clear for the second series. As to the first, it can be written as
 \begin{eqnarray*}
\nonumber&& \sum \zg_n \intx A_n(x,s)Y_n(s)\ZD s= \sum \zg_n \intx A(x,s)Y_n(s)\ZD s\\
\ZLA{eq:SeRiePerDeriYn}
&&-\sum \frac{\zg_n}{\n}\intx  K(s)\sin \n(x-s)Y_n(s)\ZD s\,.
\end{eqnarray*}
As to the second series we have the following result, used also below:

\begin{Lemma}\ZLA{Lemma:conveSERIEprostaf}{\em 
The series $ \left \{\sum (\zg_n/\n)\intx  K(s)\sin \n(x-s)Y_n(s)\ZD s \right\} $ converges in $ L^2(0,S) $.
}
\end{Lemma}
\zProof
We note that
\begin{eqnarray*}
&&\sum \frac{\zg_n}{\n}\intx  K(s)\sin \n(x-s)Y_n(s)\ZD s\\
&&=\sum \frac{\zg_n}{\n}\intx K(s)\left [ Y_n(s)-y_0\cos \n s\right ]\sin \n(x-s)\ZD s\\
&&+
y_0\sum \frac{\zg_n}{\n}\intx K(s)\sin \n(x-s)\cos \n s\ZD s\,.
 \end{eqnarray*} 
The derivative of the first series converges uniformly, thanks to 
inequality~(\ref{eq:VicinaNZAdiYn}). Using trigonometric formulas we see that the last series is

\[ 
\left ( \intx K(s)\ZD s \right )\sum \frac{\zg_n}{2\n}\sin \n x+\sum \frac{\zg_n}{2\n}\int_{-x}^x K  ( (x-r)/2)\sin \n r\ZD r
 \]
 whose derivative is $ L^2 $-convergent.\zdia
 
  We sum up: $\zsumU \zg_n\cos \n x\in W^{1,2}(0,S)  $
 and, proceeding as in previous steps, we see that the following lemma holds:
 \begin{Lemma}{\em 
 Let equality~(\ref{eq:PseKAYnFOSedipende}) hold for the sequence $ \{\zg_n\} $. 
Then, there exists $ \{\ZSI_n\}\in l^2 $ such that $ \zg_n=\ZSI_n/\n $.}
\end{Lemma}

Now we start a bootstrap argument:
\begin{Lemma}{\em 
The series $ \sum (\ZSI_n/\n) Y_n(x)$ is uniformly convergent to a $ W^{1,2}(0,S) $ function, and its derivative can be computed twice termwise.
}
\end{Lemma}
\zProof We insert $ \zg_n=\ZSI_n/\n $ in~(\ref{eq:PseKAYnFOSedipende}) and we use the representation~(\ref{eq:PrimRApprREYn}) to get
\begin{eqnarray*}
&&
\zsumU \frac{\ZSI_n}{\n}g_n(x)= \zsumU \frac{\ZSI_n}{\n^2}\intx K(s)\sin \n(x-s) Y_n(s)\ZD s\\
&&-\zsumU \frac{\ZSI_n}{\n^2}\intx \sin \n(x-s) \left [\ints A _{,11}(s,\nu) Y_n(\nu)\ZD\nu\right ]\ZD s
\end{eqnarray*}
We compute a first derivative termwise and we get  
\begin{eqnarray*}
&&\zsumU \frac{\ZSI_n}{\n}g_n'(x)=\zsumU \frac{\ZSI_n}{\n}\intx K(s)\cos \n(x-s)Y_n(s)\ZD s\\
&& -\zsumU \frac{\ZSI_n}{\n}\intx \left [\int_0^{x-\nu}
 A_{,11}(x-r,\nu)\cos \n r \ZD r
\right ] Y_n(\nu)\ZD\nu
\end{eqnarray*}
Clearly we can compute a second derivative of the last series. We prove that also the first series can be again differentiated termwise. In fact, formal differentiation gives
\[
K(x)\left (\zsumU \frac{\ZSI_n}{\n} Y_n(x)\right )-\zsumU \ZSI_n\intx K(s)\sin \n(x-s)Y_n(s)\ZD s\,.
\]
The first series is uniformly convergent and the second one is handles as in Lemma~\ref{Lemma:conveSERIEprostaf}.

The previous computations prove $ L^2 $-convergence of $ \zsumU (\ZSI_n/\n)g''_n(x) $ so that, as above:
\begin{Lemma}{\em 
If condition~(\ref{eq:PseKAYnFOSedipende}) holds then there exists a sequence $ \{\ZDE_n\}\in l^2 $ such that
\[ 
\zg_n=\frac{\ZDE_n}{\n^2}\,.
 \]
 }
\end{Lemma}
This lemma combined with~(\ref{eq:FinalePERyNcompressa}) shows that the second derivative of the series in~(\ref{eq:PseKAYnFOSedipende}) can be computed termwise, and of course it is zero. So, from~(\ref{eq:FinalePERyNcompressa}) we get
\[ 
\zsumU \ZDE_n Y_n(x)+\intx A(x,r) \left [\zsumU \ZDE_n Y_n(r)\right ]\ZD r=0\quad {\rm i.e.} \quad \zsumU \ZDE_n Y_n(x)=0\,.
 \] 
This equality, combined with~(\ref{eq:PseKAYnFOSedipende}), gives
\[ 
\sum _{n=2}^{+\ZIN}\left (1-\frac{1}{\n^2}\right ) \ZDE_n Y_n(x)=0\,.
 \]
 This argument can be repeated till we remove $ N $ first elements from the series in~(\ref{eq:PseKAYnFOSedipende}). Using $ 1- 1/\n^2\neq 0$ for $ n $ large  we conclude that the series in~(\ref{eq:PseKAYnFOSedipende}) is in fact a \emph{finite sum.} But, then all its coefficients have to be zero, since we proved linear independence of the sequence $ \{ Y_n(x)\} $.
 
 This completes the proof of Theorem~\ref{teo:Rieszbothsequences}, hence also of Theorem~\ref{teo:ControllabINTEMpot}.

 Using this result, we can now study the identification problem.
%%%%%%%%%%%%%%%%%%%%

\section{\ZLA{sec:identification}Source reconstruction}

It is now convenient to go back to the full notation $ z_n(s;t) $ and we recall that $ z_n(s;t) $ is defined for $ 0\leq s\leq t $. So, the solution of Eq.~(\ref{eq:insECONDordine}) with boundary control $ f(t)\equiv 0 $ is given by (see~(\ref{eq:RAPPRESEseriediW}))

\[ 
w(\xi,t)=\zsumU \phi_n(\xi) b_n \intt z_n(t-s;t)g(s)\ZD s\,,\quad b_n=\ZL b,\phi_n\ZR=\intp b(\xi)\phi_n(\xi)\ZD \xi\,. 
 \]
  
The output $ y(t)=w_x(0,t) $ is
\begin{equation}\ZLA{eq:SerDELLuscita}
\eta(t)=\zsumU b_n\fn \intt z_n(t-r;t)g(r)\ZD r
 \end{equation}
 We prove that we can compute the output for every $ t $, i.e.:
 \begin{Lemma}\ZLA{eq:regolarUSCITA}{\em 
The function $ t\to \eta(t) $ is continuous.
}
\end{Lemma}
\zProof

The proof is in the Appendix.

Now we choose $ g(t) $ of the special form
\begin{equation}
\ZLA{eq:GdiformaSPECIALE}
g(t)=\intt N(t-s) f(s)\ZD s\,.
\end{equation}
Here  $ f(t) $ is locally square integrable and this is achieved by taking 
\[ 
\ZSI(t)=f(t)+\intt N'(t-s) f(s)\ZD s 
 \]
 (use $ N(0)=1 $.)
Let us fix any $ T\geq T_0 $ ($ T_0 $ is specified in Theorem~\ref{teo:ControllabINTEMpot}) and let us note that
\begin{eqnarray*}
&&\eta(T)=\left \ZL \zsumU \phi_n(x) \fn\intT z_n(T-r;T)\intr N(r-s) f(s)\ZD s\,\ZD r,b_n(x)\right \ZR\\
&&= \left \ZL \zsumU \phi_n(x)\intT f(T-s) \left [
\fn\ints N(s-r)z_n(r;T)\ZD r
\right ]\ZD s,b(x)\right \ZR\ 
\end{eqnarray*}
(the crochet denotes $ L^2(0,\pi) $-inner product.)
We recall $ z_n(r;T)=z_n(r) $ and we compare this formula with~(\ref{eq:seriediW}).
We see that the left side of the crochet is nothing else then $ w(T) $ when $ f(t) $ is the \emph{boundary control.} Hence, for every $ k $ there exists a suitable $ f_k(t)\in L^2(0,T) $, i.e. a suitable $ \ZSI_k(t)\in L^2(0,T) $ such that the corresponding output $ \eta_{(k)}(t) $ is such that
\[ 
\eta_{(k)}(T)=\ZL \phi_k,b\ZR=b_k\,,
 \]
 the $ k-th $ Fourier coefficient of $ b(\xi) $ respect to the orthonormal basis $ \{\phi_k(\xi)\} $ (or, if we whish, to any prescribed orthonormal basis, thanks to
  Theorem~\ref{teo:ControllabINTEMpot}.) So, $ b(\xi) $ is given by
 \[ 
 b(\xi)=\zsumU \eta_{(k)}(T) \phi_k(\xi)\,.
  \] 
 \begin{Remark}
{ \em
At first sight it might seem that the previous method is similar to the one developed in~\cite{Yamam}. In fact, it is not the same idea and, as has to be expected due to the time varying coefficient $ P(t) $, it is less efficient. In fact, the algorithm in~\cite{Yamam} does not really use $ \ZSI(t) $ which has to be known, smooth and with $ \ZSI(0)=1 $, but fixed once and for all. Using this fixed input the output $ y(t) $  is measured for every $ t\in[0,T] $ and then a sequence of algorithms applied to this simple observation produce the Fourier coefficients $ b_k $. The algorithm in~\cite{Yamam} has been extended to systems with memory, but with constant $ P(t) $, in~\cite{PandIDENT} but it seems that this algorithm cannot be used to study the time varying traction. For this reason we proposed to use the full map $   L^2(0,T)\ni\ZSI\ZCD \mapsto y(T) $.\zdia
}
\end{Remark}

\section{Appendix: The proof of Lemmas~\ref{eq:Lemma:regol-soluzione},   \ref{lemmaVicinaNZAdiYn} and~\ref{eq:regolarUSCITA}}

In this appendix it is important to recall that every transformation in Section~\ref{sect:trabNSf} is on a fixed interval $ [0,T] $ and {\em does depend on $ T $.\/} Even the transformation $ H(t)=H(t;T) $ since it depends on $ Q(t)=P(T-t) $.  Now we considere these transformation for every fixed $ T $ in an interval $ [0,K] $ so that it will be $ 0\leq t\leq T\leq K $ and the pair $ (t,T) $ will belong to the triangle 
\[ 
 \tilde \triangle= \{(t,T)\,:\quad 0\leq t\leq T\leq K\}\,.
  \]
   We have that 
 \[ 
 L(t)=L(t;T)=\int_0^{t} \sqrt{P(T -s)}\ZD s
  \]
 is bounded on the triangle $ \tilde \triangle $ and so also is $ L'(t;T) $, uniformly for $ 0\leq T\leq K $. 
  
  As in Section~\ref{sect:trabNSf}, let, for every $ T\in[0,K] $,
  \[ 
  S=L(T)=L(T;T)
   \]  
 Note that $S= L(T;T) $ so that $ 0\leq x= L(t;T)\leq S\leq \Xi=L(K,K) $
 and that the triangle $ \tilde \triangle $ is transformed to
 \[ 
  \Delta=\{(x,S)\,:\quad 0\leq x\leq S\leq \Xi\}\,.
   \]
 
  So,  $ M(t;X) $ and $ M'(t;X) $ and $ A(x,r;X) $ and their derivatives  are bounded of $ \triangle $.
  We recall that $ y_0 $ and $ y_1 $ do depend on $ S $. So we shall write
   \[
    g_n(x)=g_n(x;S) =y_0(S)\cos\n x+\frac{1}{\fn}y_1(S)\sin\n x
    \]
    and that $ y_0(S) $ and $y_1(S)$ are bounded on $ [0,\Xi] $.

\subparagraph{Proof of Lemma~\ref{lemmaVicinaNZAdiYn}.}
The first inequality in~(\ref{eq:VicinaNZAdiYn}) follows from Gronwall inequality applied to~(\ref{eq:PrimRApprREYn}). So, we must prove
\begin{equation}\ZLA{eq:DiseqLEMMAinAppendice}
|Y_n(x;S)-g_n(x;S)|\leq \frac{M}{n } 
\end{equation}
 (the constant $ M $ depends on $ \Xi $.)
  
  In order to get inequality~(\ref{eq:DiseqLEMMAinAppendice}) we integrate by parts once the last integral in~(\ref{eq:PrimRApprREYn}) and we add and subtract $ g_n(x;S) $ to the first integral. Then we integrate by parts so to get a factor $ 1/\n $. We obtain
  
 \begin{eqnarray*}
\nonumber&& Y_n(x)-g_n(x;S)=
\\&& -\intx A_n(x,r;S)\left [Y_n(r;S)-g_n(r;S)\right ]\ZD s
 \\&&-\intx A_n(x,r;S)g_n(s;S)\ZD r+\frac{1}{\n^2}M_n(x;S)
\end{eqnarray*}
where $ M_n(x;S) $ is a function that we don't need to write down explicitly, but such that $ |M_n(x,S)|<M $ for every $ n $ and $ (x,S)\in\triangle $.

The require inequality follows from here, integrating by parts the last integral and using Gronwall inequality.\zdia

  \subparagraph{The proofs of Lemmas~\ref{eq:Lemma:regol-soluzione} and ~\ref{eq:regolarUSCITA}.} 
  
  The idea of the proofs of these lemmas is similar. In order to prove Lemma ~\ref{eq:Lemma:regol-soluzione}
we have to consider   the first series in~(\ref{eq:RAPPRESEseriediW})
while for Lemma~\ref{eq:regolarUSCITA} we have to consider the series
 in~(\ref{eq:SerDELLuscita}). The common feature is an integral of the general form
 
\begin{eqnarray*}
&& \intt V(t,r )z_n(r;T)\ZD r \\
&&
= \int_0^{L(t;T)}
   V(t;M(s;S) )  e^{N'(0)M(s;S)} P^{-1/4}(M(s;S)-T)\tilde\zeta_n(M(s;S);T)\ZD s\,.
 \end{eqnarray*}
Now we use  
\[ 
L(M(x;S);T)=x\,,\quad Y_{  n} (s;S)=\tilde \zeta_n(M(s;S);T)
 \]
and finally we get    an  integral of the form

\CM{``an'' era scritto ``and''}
\[ 
\int_0^x V_1(x,s;X)Y_{ n} (s;X)\ZD s\,.
 \]
The integral in Lemma~\ref{eq:Lemma:regol-soluzione} is obtained when
$ V(t,r;T)=N(t-r) $ and then $ V_1(x,s;X)=C(x,s;X) $. The integral in~(\ref{eq:SerDELLuscita})  is obtained with $ V(t,r ) =g(t-r)$. Note that $ C(x,x;X)\neq 0 $ while $ g(0)=0 $ and so in the proof of Lemma~\ref{eq:regolarUSCITA} we shall have $ V_1(s,s;X)=0 $.

We present now the proof of Lemma~\ref{eq:Lemma:regol-soluzione} and
 we leave the similar proof (based on inequality~\ref{eq:DISEQperZ})  of Lemma~\ref{eq:regolarUSCITA} to the reader.

 We transform the variable $ t \in [0,T]\subseteq [0,K]$ in the series~(\ref{eq:RAPPRESEseriediW}) 
to the variable $ x=M(t)\in [0,S]\subseteq [0,\Xi] $  as shown above  and we consider the $ L^2(0,\pi) $
norm of the resulting series, for each $ S\in[0,\Xi] $. Using the fact that $ \{\phi_n(\xi)\} $ is orthonormal in $ L^2(0,\pi) $ and definition~(\ref{eq:definCxs}), the square of the norm can be written as
\begin{eqnarray*}
&&
\sum _{n=1}^{+\ZIN}\left |\int_0^x \Lambda(x,\nu;S)\left (\fn \int_0^{\nu} C(\nu,r;S) Y_n(r;S)\ZD r\right )\ZD\nu 
\right |^2\\
&&=
\sum _{n=1}^{+\ZIN}\left |\int_0^x \Lambda(x,\nu;S)\left (C(\nu,\nu;S)Z_n(\nu;S)\right )\ZD\nu 
\right |^2
 \\
&&
\Lambda(x,\nu;S)=\left ( f\left (M(x;S )-M(\nu;S )\right )M'(\nu;S )\right )
 \end{eqnarray*}
Using~(\ref{eq:DISEQperZ}) we see that it is sufficient that we note
\[ 
\sum _{n=1}^{+\ZIN}\left |\int_0^x \Lambda(x,\nu;S)C(\nu,\nu;S)\sin ns\right |^2\leq \int_0^x \left |\Lambda(x,\nu;S)C(\nu,\nu;S)\right |^2\ZD s\,.
 \]

The proof is now finished since the right hand side is bounded for $ 0\leq s\leq S\leq\Xi $.\zdia

\end{document}